\documentclass[a4paper,12pt]{article}
\usepackage{amsmath}
\usepackage{amssymb}
\usepackage{amsfonts}
\usepackage{latexsym}
\usepackage{enumerate}
\usepackage{stmaryrd}

\newtheorem{llemma}{Lemma}[section]

\newtheorem{exmp}[llemma]{Example}
\newtheorem{thm}[llemma]{Theorem}

\newtheorem{defn}[llemma]{Definition}
\newtheorem{rem}[llemma]{Remark}

\newtheorem{key}[llemma]{Keyword}

\newcommand{\proof}{\emph{Proof. }}

\begin{document}
\title{\textbf{A new approach to bipolar soft sets and its applications}}
\author{Faruk Karaaslan$^a$\thanks{corresponding author}\,\,and  Serkan Karata\c{s}$^b$ }
\date{\footnotesize{$^a$Department of Mathematics, Faculty of Sciences, Karatekin
University, 18100 \c{C}ank{\i}r{\i}, Turkey,
fkaraaslan@karatekin.edu.tr\\
$^b$Department of Mathematics, Faculty of Arts and Sciences, Ordu
University, 52200 Ordu, Turkey serkankaratas@odu.edu.tr}}

\maketitle

%***********************************************************************************
% ABSTRACT
\begin{abstract}
Molodtsov \cite{molodtsov-1999} proposed the concept of soft set theory in 1999, which can be used as a  mathematical tool  for dealing with problems that contain uncertainty. Sabir and Naz \cite{shabir-2013} defined notion of bipolar  soft set  in 2013. In this paper, we redefine concept of bipolar soft set and bipolar soft set operations  as more functional than former. Also  we study on  their basic properties and we present a decision making method with application.

\begin{key}
Soft set, bipolar soft set, decision making, bipolar soft set operations
\end{key}
\end{abstract}

%***********************************************************************************
%% section 1: introduction

\section{Introduction}
In the real life, problems in economy, engineering, environmental science and
social science and many fields involve data that contain uncertainties.
These problems may not be successfully modeled by existing methods in
classical mathematics. There are some well known mathematical theories
for dealing with uncertainties such as; fuzzy set theory \cite{zadeh-1965},
intuitionistic fuzzy set theory \cite{atanassov-1986}, rough set \cite{pawlak-1982}
and so on. But all of these theories have their own difficulties which are pointed
out in \cite{molodtsov-1999}.  To cope with these difficulties, Molodtsov proposed
the concept of soft set as a new mathematical tool for dealing with uncertainties.
Then Maji et al. \cite{maji-2003},  equality of two soft sets, subset and super
set of a soft set, complement of a soft set, null soft set, and absolute soft
set with examples. Moreover, they defined soft binary operations like AND, OR
and  the operations of union, intersection. In 2009, Ali et al.\cite{ali-09}
gave some new notions such as the restricted intersection, the restricted
union, the restricted difference and the extended  intersection of two soft sets.
In 2010, to make soft set operations more functional, \c{C}a\u{g}man and Engino\u{g}lu
 \cite{cagman-2009}  redefine the operations of Molodtsov's soft sets. Also
they defined products of soft sets and uni–int decision function and
constructed an uni-int decision making method which selects a set of optimum elements
from the alternatives. Feng et al. \cite{feng-2010a,feng-2011,feng-2010b} combined
 soft sets fuzzy sets and rough sets and gave application in decision making.
In 2014, \c{C}a\u{g}man \cite{cagman-2014}  made contributions to the theory of
soft sets to fill gaps of former definition and operations.
Recently, the properties and applications on the soft set theory have been studied increasingly \cite{aktas-2007,cagman-2010,kharal-2010,qin-10,kong-2014,maju-8,min-12}.
By embedding the ideas of fuzzy sets, intuitionistic fuzzy sets and interval-valued
fuzzy sets, many interesting applications of soft set theory have been expanded
\cite{cagman-2011a,cagman-2011b,cagman-2013,deli-2012a,deli-2012b,jiag-2010,jiag-2011,kar-2014,maji-2001,som-2006,sut-2012,xiao-2009,yang-2009,zou-2008b}.
In 2013, Sabir and Naz \cite{shabir-2013} defined bipolar soft sets and basic
operations of union, intersection and complementation  for bipolar soft sets.
They gave examples of bipolar soft sets and an application of bipolar soft sets
in a decision making problem. In 2013, Aslam et al. \cite{aslam-2013} introduced the
notion of bipolar fuzzy soft set and study fundamental properties.
Also, Naz and Sabir \cite{naz-2014}  studied on algebraic structure of bipolar fuzzy soft sets.

In this paper,  based on \c{C}a\u{g}man \cite{cagman-2014}, we redefine the
notion of bipolar soft sets and bipolar soft set operations
as more functional. Since it is difficult  to express 'not' of an element as mathematical,
to cope with this case we use to define a bijective  function.
Also we give decision making method for decision making problem
 involving the bipolar soft set and present an example related to this method.

%***********************************************************************************
%% Section 2: preliminaries
\section{Preliminaries}
In this section, we will remind the definition of soft set which is
defined by Molodtsov \cite{molodtsov-1999} and present contributions
to soft set which is made by \c{C}a\u{g}man \cite{cagman-2014}.
Moreover, we will give some properties of these topics. Throughout
this paper, we will denote initial universe, set of parameters and
power set of $U$ by $U$, $E$ and $P(U)$, respectively.

%-----------------------------------------------------------------------------------
\begin{defn}\cite{molodtsov-1999}
Let consider a nonempty set $A$, $A\subseteq E$. A pair $(F,A)$ is
called a soft set over $U$, where $F$ is a mapping given by $F:A\to
P(U)$.
\end{defn}

In this paper, we will use following definition which was defined by \c{C}a\u{g}man
\cite{cagman-2014} and we will do our work based on this definition.

%-----------------------------------------------------------------------------------
\begin{defn}\cite{cagman-2014}
A soft set $F$ over $U$ is a set valued function from $E$ to $\mathcal{P}(U)$.  It can
be written a set of ordered pairs
$$
F=\big\{(e,F(e)):e\in E\big\}.
$$
Note that if $F(e)=\emptyset$, then the element $(e, F(e))$ is not appeared in $F$.
Set of all soft sets over $U$ is denoted by $\mathbb{S}$.
\end{defn}

%-----------------------------------------------------------------------------------
\begin{defn}\cite{cagman-2014} Let $F,G\in\mathbb{S}$. Then,
\begin{enumerate}[\it i.]
\item If $F(e)=\emptyset$ for all $e\in E$, $F$ is said to be a null soft set, denoted
by $\Phi$.
\item If $F(e)=U$ for all $e\in E$, $F$ is said to be absolute soft set, denoted by
$\hat{U}$.
\item $F$ is soft subset of $G$, denoted by $F\tilde\subseteq G$, if $F(e)\subseteq G(e)$
for all $e\in E$.
\item $F=G$, if $F\tilde\subseteq G$ and $G\tilde\subseteq F$.
\item Soft union of $F$ and $G$, denoted by $F\tilde\cup G$, is a soft set over $U$ and
defined by $F\tilde\cup G:E\to P(U)$ such that $(F\tilde\cup
G)(e)=F(e)\cup G(e)$ for all $e\in E$.
\item Soft intersection of $F$ and $G$, denoted by $F\tilde\cap G$, is a soft set over
$U$ and defined by $F\tilde\cap G:E\to P(U)$ such that $(F\tilde\cap
G)(e)=F(e)\cap G(e)$ for all $e\in E$.
\item Soft complement of $F$ is denoted by $F^{\tilde c}$ and defined by
$F^{\tilde c}:E\to P(U)$ such that $F^{\tilde c}(e)=U\setminus F(e)$
for all $e\in E$.
\end{enumerate}
\end{defn}

%***********************************************************************************
%% section: Bipolar soft sets
\section{Bipolar soft sets}
In this section, we will redefine bipolar soft set and introduce some basic properties.

%\begin{defn} Let $E$ be a parameter set and $E_2$ be NOT set  of
%$E_1$. Then, $E_1\cup E_2$  is called parameter set extending with
%parameter sets $E_1$ and $E_2$ and denoted by
%$\mathcal{E}_{E_1,E_2}$.
%\end{defn}
%\begin{defn} Let $E_1$ and $E_2$ be two parameter (attribute) sets
%that have the same number of elements and $f$ be a bijective
%function from $E_1$ to  $E_2$ such that $f(e)=\neg e$, for all $e\in
%E_1$. Then, extending parameter set, denoted by $\mathcal{E}$, is
%defined by as follow
%$$
%\mathcal{E}=E_1\cup E_2
%$$
%\end{defn}
%\begin{defn}Let $\mathcal{E}_{E_1,E_2}$ be parameter set  extending  with parameter sets $E_1$ and $E_2$. If $F:E_1\to P(U)$ and $G:E_2\to
%P(U)$ are two mappings such that $F(e)\cap G(f(e))=\emptyset$, then
%triple $(F,G,E)$ is called bipolar soft set, where  $f:E_1\to E_2$
%is a bijective function which defined by $f(e_i)=e_j$ that $e_j=\neg
%e_i$, for all $e_i \in E_1$. Set of all bipolar soft sets over $U$
%is denoted by $\mathbb{BS}$. We can represent a bipolar soft set
%$(F,G,E)$ as following form:
%$$
%(F,G,E)=\Big\{\big\langle(e,F(e)),(f(e),G(f(e)))\big\rangle:e\in E_1
%\,\,\,\text{and}\,\,\,F(e)\cap G(f(e))=\emptyset\Big\}.
%$$
%If $F(e)=\emptyset$ and $G(f(e))=\emptyset$ for $e\in E_1$, then
%$\big\langle(e,\emptyset),(f(e),\emptyset)\big\rangle$ is no written
%in the bipolar soft set $(F,G,E)$.
%\end{defn}
%-----------------------------------------------------------------------------------
\begin{defn} Let $E$ be a parameter set and  $E_1$ and $E_2$ be two non empty subsets of $E$
such that $E_1\cup E_2=E$ and $E_1\cap E_2=\emptyset$. If $F:E_1\to
P(U)$ and  $G:E_2\to P(U)$ are two mappings such that $F(e)\cap
G(f(e))=\emptyset$, then triple $(F,G,E)$ is called bipolar soft
set, where  $f:E_1\to E_2$ is a bijective function. Set of all
bipolar soft sets over $U$ is denoted by $\mathbb{BS}$. We can
represent a bipolar soft set $(F,G,E)$ as following form:
$$
(F,G,E)=\Big\{\big\langle(e,F(e)),(f(e),G(f(e)))\big\rangle:e\in E_1
\,\,\,\text{and}\,\,\,F(e)\cap G(f(e))=\emptyset\Big\}.
$$
If $F(e)=\emptyset$ and $G(f(e))=\emptyset$ for $e\in E_1$, then
$\big\langle(e,\emptyset),(f(e),\emptyset)\big\rangle$ is no written in the bipolar
soft set $(F,G,E)$.
\end{defn}
%% There is no mathematical provision that $\neg e$.

%-----------------------------------------------------------------------------------
\begin{rem}
According to above form of the bipolar soft set $(F,G,E)$, we can construct a tabular
form every bipolar soft set for a finite universal set $U$ and finite parameter set
$E$. Therefore,  let $U=\{u_1,u_2,\ldots,u_m\}$ and $E=\{e_1,e_2,\ldots,e_n\}$. Then,
\begin{center}
\begin{tabular}{c|cccc}
$(F,G,E)$ & $(F(e_1),G(f(e_1)))$ & $(F(e_2),G(f(e_2)))$ & $\ldots$ & $(F(e_n),G(f(e_n)))$ \\
  \hline
$u_1$ & $(a_{11},b_{11})$ & $(a_{12},b_{12})$ & $\ldots$ & $(a_{1n},b_{1n})$ \\
$u_2$ & $(a_{21},b_{21})$ & $(a_{22},b_{22})$ & $\ldots$ & $(a_{2n},b_{2n})$ \\
$\vdots$ & $\vdots$ & $\vdots$ & $\ddots$ & $\vdots$\\
$u_m$ & $(a_{m1},b_{m1})$ & $(a_{m2},b_{m2})$ & $\ldots$ & $(a_{mn},b_{mn})$ \\
\end{tabular}
\end{center}
where,
$$
a_{ij}=
\begin{cases}
1,& u_i\in F(e_j)\\
0,& u_i\notin F(e_j)
\end{cases}
\,\,\,\text{and}\,\,\,
b_{ij}=
\begin{cases}
1,& u_i\in G(f(e_j))\\
0,& u_i\notin G(f(e_j)).
\end{cases}
$$
Note that, $a_{ij}$ and $b_{ij}$ must not be $1$ in same time. So, there are three
cases for $(a_{ij},b_{ij})$: $(1,0)$, $(0,1)$ or $(0,0)$.
\end{rem}
%-----------------------------------------------------------------
\begin{exmp}\label{e-bipolarsoftset}
Let $U=\{u_l,u_2,u_3,u_4,u_5,u_6,u_7,u_8\}$ be the universe which are eight houses
and $E=\{e_1,e_2,e_3,e_4,e_5,e_6,e_7,e_8\}$ be the set of parameters. Here, $e_i$
$(i=1,2,3,4,5,6,7,8)$ stand for the parameters  ``\emph{large}'', ``\emph{small}'',
``\emph{modern}'', ``\emph{standard}'', ``\emph{cheap}'', ``\emph{expensive}''
``\emph{with parking}'',  and ``\emph{no parking area}''  respectively.  Therefore,
we can chose $E_1$ and $E_2$ sets as $E_1=\{e_1,e_3,e_5,e_7\}$ and $E_2=\{e_2,e_4,e_6,e_8\}$.
Now, we define the bijective function $f$ as $f(e_i)=\neg e_i$ $(i=1,3,5,7)$. Here,
the notion $\neg e_i$ means that ``\emph{not} $e_i$'' for all $i=1,3,5,7$. Thus,
we have following results.
\begin{eqnarray*}
& & f(e_1) = \neg e_1=e_2 \\
& & f(e_3) = \neg e_3=e_4 \\
& & f(e_5) = \neg e_5=e_6 \\
& & f(e_7) = \neg e_7=e_8
\end{eqnarray*}
So, we can describe following t bipolar soft sets $(F_1,G_1,E)$ and $(F_2,G_2,E)$
to buy a house.
\begin{eqnarray*}
(F_1,G_1,E) & = &
\Big\{
\big\langle(e_1,\{u_1,u_3,u_4\}),(e_2,\{u_2,u_6\})\big\rangle,\\
&&\big\langle(e_3,\{u_2,u_5,u_7\}),(e_4,\{u_1,u_3,u_8\})\big\rangle,\\
&&\big\langle(e_5,\{u_3,u_4\}),(e_6,\{u_1,u_2,u_5,u_8\})\big\rangle,\\
&&\big\langle(e_7,\{u_5,u_6,u_7,u_8\}),(e_8,\{u_2,u_3\})\big\rangle
\Big\}
\end{eqnarray*}
\begin{eqnarray*}
(F_2,G_2,E) & = &
\Big\{
\big\langle(e_1,\{u_1,u_2,u_4\}),(e_2,\{u_3,u_5,u_6,u_7\})\big\rangle,\\
&&\big\langle(e_3,\{u_2,u_5\}),(e_4,\{u_1,u_3,u_4,u_8\})\big\rangle,\\
&&\big\langle(e_5,\{u_1,u_3,u_4\}),(e_6,\{u_2,u_5,u_7,u_8\})\big\rangle,\\
&&\big\langle(e_7,\{u_5\}),(e_8,\{u_2,u_3,u_4\})\big\rangle
\Big\}
\end{eqnarray*}
\end{exmp}
Tabular representations of bipolar soft sets $(F_1,G_1,E)$ and $(F_2,G_2,E)$ are in table ....
\scriptsize
\begin{center}
\begin{tabular}{c|cccc}
$(F_1,G_1,E)$ & $(F_1(e_1),G_1(f(e_1)))$ & $(F_1(e_3),G_1(f(e_3)))$ & $(F_1(e_5),G_1(f(e_5)))$ & $(F_1(e_7),G(f_1(e_7)))$ \\
  \hline
$u_1$ & $(1,0)$ & $(0,1)$ & $(0,1)$ & $(0,0)$ \\
$u_2$ & $(0,1)$ & $(1,0)$ & $(0,1)$ & $(0,1)$ \\
$u_3$ & $(1,0)$ & $(0,1)$ & $(1,0)$ & $(0,1)$\\
$u_4$ & $(1,0)$ & $(0,0)$ & $(1,0)$ & $(0,0)$ \\
$u_5$ & $(0,0)$ & $(1,0)$ & $(0,1)$ & $(1,0)$ \\
$u_6$ & $(0,1)$ & $(0,0)$ & $(0,0)$ & $(1,0)$ \\
$u_7$ & $(0,0)$ & $(1,0)$ & $(0,0)$ & $(1,0)$\\
$u_8$ & $(0,0)$ & $(0,1)$ & $(0,1)$ & $(1,0)$ \\
\end{tabular}
\end{center}
\normalsize and

\scriptsize
\begin{center}
\begin{tabular}{c|cccc}
$(F_2,G_2,E)$ & $(F_2(e_1),G_2(f(e_1)))$ & $(F_2(e_3),G_2(f(e_3)))$ & $(F_2(e_5),G_2(f(e_5)))$ & $(F_2(e_7),G(f_2(e_7)))$ \\
  \hline
$u_1$ & $(1,0)$ & $(0,1)$ & $(1,0)$ & $(0,0)$ \\
$u_2$ & $(1,0)$ & $(1,0)$ & $(0,1)$ & $(0,1)$ \\
$u_3$ & $(0,1)$ & $(0,1)$ & $(1,0)$ & $(0,1)$\\
$u_4$ & $(1,0)$ & $(0,1)$ & $(1,0)$ & $(0,1)$ \\
$u_5$ & $(0,1)$ & $(0,0)$ & $(0,1)$ & $(1,0)$ \\
$u_6$ & $(0,1)$ & $(0,0)$ & $(0,0)$ & $(0,0)$ \\
$u_7$ & $(0,1)$ & $(0,0)$ & $(0,1)$ & $(0,0)$\\
$u_8$ & $(0,0)$ & $(0,1)$ & $(0,1)$ & $(0,0)$ \\
\end{tabular}
\end{center}
\normalsize
%-----------------------------------------------------------------------------------
\begin{defn}\label{d-subset}
Let $(F_1,G_1,E),(F_2,G_2,E)\in\mathbb{BS}$. Then,  $(F_1,G_1,E)$ is bipolar soft
subset of $(F_2,G_2,E)$, denoted by $(F_1,G_1,E)\sqsubseteq(F_2,G_2,E)$, if
$F_1(e)\subseteq F_2(e)$  and $G_2(f(e))\subseteq G_1(f(e))$ for all $e\in E_1$.
\end{defn}
%------------------------------------------------------------------------------------
\begin{defn}
Let $(F_1,G_1,E),(F_2,G_2,E)\in\mathbb{BS}$. Then,  $(F_1,G_1,E)$ and $(F_2,G_2,E)$
equal, denoted by $(F_1,G_1,E)=(F_2,G_2,E)$, if $(F_1,G_1,E)\sqsubseteq(F_2,G_2,E)$
and $(F_2,G_2,E)\sqsubseteq(F_1,G_1,E)$.
\end{defn}

%-----------------------------------------------------------------------------------
\begin{defn} \label{d-union}
Let $(F_1,G_1,E),(F_2,G_2,E)\in\mathbb{BS}$. Then,  bipolar soft
union of $(F_1,G_1,E)$ and $(F_2,G_2,E)$ is, denoted by $(F_1,G_1,E)\sqcup(F_2,G_2,E)$,
defined by $(F_1\sqcup F_2)(e)=F_1(e)\cup F_2(e)$  and
$(G_1\sqcup G_2)(f(e))=G_1(f(e))\cap G_2(f(e))$ for all $e\in E_1$.
\end{defn}

%-----------------------------------------------------------------------------------
\begin{defn} \label{d-intersection}
Let $(F_1,G_1,E),(F_2,G_2,E)\in\mathbb{BS}$. Then,  bipolar soft intersection of
$(F_1,G_1,E)$ and $(F_2,G_2,E)$ is, denoted by $(F_1,G_1,E)\sqcap(F_2,G_2,E)$,
defined by $(F_1\sqcap F_2)(e)=F_1(e)\cap F_2(e)$  and
$(G_1\sqcap G_2)(f(e))=G_1(f(e))\cup G_2(f(e))$ for all $e\in E_1$.
\end{defn}

%-----------------------------------------------------------------------------------
\begin{defn} \label{d-null}
Let $(F,G,E)\in\mathbb{BS}$. For all $e\in E_1$, $F(e)=\emptyset$ and $G(f(e))=U$,
then $(F,G,E)$ is called  null bipolar soft set and denoted by $(\Phi,\tilde U,E)$.
\end{defn}

%-----------------------------------------------------------------------------------
\begin{defn} \label{d-absolute}
Let $(F,G,E)\in\mathbb{BS}$. For all $e\in E_1$, $F(e)=U$ and $G(f(e))=\emptyset$,
then $(F,G,E)$ is called  absolute bipolar soft set and denoted by $(\tilde U,\Phi,E)$.
\end{defn}

%-----------------------------------------------------------------------------------
\begin{defn} \label{d-complement}
Let $(F,G,E)\in\mathbb{BS}$. Then, complement of $(F,G,E)$, denoted by $(F,G,E)^{\tilde c}$,
is a bipolar soft set over $U$ such that $(F,G,E)^{\tilde c}=(H,K,E)$, where
$H(e)=G(f(e))$ and $K(f(e))=F(e)$ for all $e\in E_1$.
\end{defn}

%-----------------------------------------------------------------------------------
\begin{exmp}\label{exsubintuni}
Let $U=\{u_l,u_2,u_3,u_4,u_5,u_6,u_7,u_8\}$ be an initial universe
and
 $E=\{e_1,e_2,e_3,e_4,e_5,e_6,e_7,e_8\}$ be a set of parameters. If we
 chose $E_1$ and $E_2$ sets as $E_1=\{e_1,e_3,e_5,e_7\}$ and $E_2=\{e_2,e_4,e_6,e_8\}$.
Now, we define the bijective function $f$ as $f(e_i)=\neg e_i$ $(i=1,3,5,7)$. Here,
the notion $\neg e_i$ means that ``\emph{not} $e_i$'' for all $i=1,3,5,7$. Thus,
we have following results.
\begin{eqnarray*}
& & f(e_1) = \neg e_1=e_2 \\
& & f(e_3) = \neg e_3=e_4 \\
& & f(e_5) = \neg e_5=e_6 \\
& & f(e_7) = \neg e_7=e_8
\end{eqnarray*}
So, we can describe following the bipolar soft sets $(F_1,G_1,E)$, $(F_2,G_2,E)$ and $(F_3,G_3,E)$.
\begin{eqnarray*}
(F_1,G_1,E) & = &
\Big\{
\big\langle(e_1,\{u_1,u_3,u_4\}),(e_2,\{u_2,u_6\})\big\rangle,\\
&&\big\langle(e_3,\{u_2,u_5,u_7\}),(e_4,\{u_1,u_3,u_8\})\big\rangle,\\
&&\big\langle(e_5,\{u_3,u_4\}),(e_6,\{u_1,u_2,u_5,u_8\})\big\rangle,\\
&&\big\langle(e_7,\{u_5,u_6,u_7,u_8\}),(e_8,\{u_2,u_3\})\big\rangle
\Big\}
\end{eqnarray*}
\begin{eqnarray*}
(F_2,G_2,E) & = &
\Big\{
\big\langle(e_1,\{u_1,u_2,u_4\}),(e_2,\{u_3,u_5,u_6,u_7\})\big\rangle,\\
&&\big\langle(e_3,\{u_2,u_5\}),(e_4,\{u_1,u_3,u_4,u_8\})\big\rangle,\\
&&\big\langle(e_5,\{u_1,u_3,u_4\}),(e_6,\{u_2,u_5,u_7,u_8\})\big\rangle,\\
&&\big\langle(e_7,\{u_5\}),(e_8,\{u_2,u_3,u_4\})\big\rangle
\Big\}
\end{eqnarray*}
\begin{eqnarray*}
(F_3,G_3,E) & = &
\Big\{
\big\langle(e_1,\{u_1,u_4\}),(e_2,\{u_2,u_3,u_5,u_6\})\big\rangle,\\
&&\big\langle(e_3,\{u_2,u_5\}),(e_4,\{u_1,u_3,u_4,u_8\})\big\rangle,\\
&&\big\langle(e_5,\{u_3\}),(e_6,\{u_1,u_2,u_5,u_7,u_8\})\big\rangle,\\
&&\big\langle(e_7,\{u_5,u_8\}),(e_8,\{u_2,u_3,u_4,u_7\})\big\rangle
\Big\}
\end{eqnarray*}
\end{exmp}
Note that, for all $e\in E_1$,\\
since $F_3(e)\subseteq F_1(e)$ and $G_3(f(e))\supseteq G_1(f(e))$,
$(F_3,G_3,E)\sqsubseteq (F_1,G_1,E)$ \\

\begin{eqnarray*}
(F_1,G_1,E)\sqcup (F_2,G_2,E)&=&\{(e_1,\{u_1,u_4\},\{u_2,u_3,u_5,u_6,u_7\}),\\
&&
(e_3,\{u_2,u_5\},\{u_1,u_3,u_4,u_8\}),
(e_5,\{u_3,u_4\},\\
&&\{u_1,u_2,u_5,u_7,u_8\}), (e_7,\{u_5\},\{u_2,u_3,u_4\})\}
\end{eqnarray*}
and \\
\begin{eqnarray*}
(F_1,G_1,E)\sqcap (F_2,G_2,E)&=&\{(e_1,\{u_1,u_2,u_3,u_4\},\{u_6\}),\\
&&(e_3,\{u_2,u_5,u_7\},\{u_1,u_3,u_8\}),
(e_5,\{u_1,u_3,u_4\},\\
&&\{u_2,u_5,u_8\}), (e_7,\{u_5,u_6,u_7,u_8\},\{u_2,u_3\})\}
\end{eqnarray*}

%-----------------------------------------------------------------------------------
\begin{thm}
Let $(F,G,E),(F_1,G_1,E),(F_2,G_2,E)\in\mathbb{BS}$. Then,
\begin{enumerate}[\it i.]
\item $(F,G,E)\sqsubseteq(\tilde U,\Phi,E)$
\item $(\Phi,\tilde U,E)\sqsubseteq(F,G,E)$
\item $(F,G,E)\sqsubseteq(F,G,E)$
\item If $(F,G,E)\sqsubseteq(F_1,G_1,E)$ and $(F_1,G_1,E)\sqsubseteq(F_2,G_2,E)$, then
$(F,G,E)\sqsubseteq(F_2,G_2,E)$
\end{enumerate}
\end{thm}
\proof It is clear from Definition \ref{d-subset}.

%-----------------------------------------------------------------------------------
\begin{thm}
Let $(F,G,E),(F_1,G_1,E),(F_2,G_2,E)\in\mathbb{BS}$. Then,
\begin{enumerate}[\it i.]
\item $(F,G,E)\sqcup(F,G,E)=(F,G,E)$
\item $(F,G,E)\sqcup(\Phi,\tilde U,E)=(F,G,E)$
\item $(F,G,E)\sqcup(\tilde U,\Phi,E)=(\tilde U,\Phi,E)$
\item $(F,G,E)\sqcup(F,G,E)^{\tilde c}=(\tilde U,\Phi,E)$
\item $(F_1,G_1,E)\sqcup(F_2,G_2,E)=(F_2,G_2,E)\sqcup(F_1,G_1,E)$
\item $(F,G,E)\sqcup\big[(F_1,G_1,E)\sqcup(F_2,G_2,E)\big]=
\big[(F,G,E)\sqcup(F_1,G_1,E)\big]\sqcup(F_2,G_2,E)$
\end{enumerate}
\end{thm}
\proof It can be proved by Definition \ref{d-union}.

%-----------------------------------------------------------------------------------
\begin{thm}
Let $(F,G,E),(F_1,G_1,E),(F_2,G_2,E)\in\mathbb{BS}$. Then,
\begin{enumerate}[\it i.]
\item $(F,G,E)\sqcap(F,G,E)=(F,G,E)$
\item $(F,G,E)\sqcap(\Phi,\tilde U,E)=(\Phi,\tilde U,E)$
\item $(F,G,E)\sqcap(\tilde U,\Phi,E)=(F,G,E)$
\item $(F,G,E)\sqcap(F,G,E)^{\tilde c}=(\Phi,\tilde U,E)$
\item $(F_1,G_1,E)\sqcap(F_2,G_2,E)=(F_2,G_2,E)\sqcap(F_1,G_1,E)$
\item $(F,G,E)\sqcap\big[(F_1,G_1,E)\sqcap(F_2,G_2,E)\big]=
\big[(F,G,E)\sqcap(F_1,G_1,E)\big]\sqcap(F_2,G_2,E)$
\end{enumerate}
\end{thm}
\proof It can be proved by Definition \ref{d-intersection}.

%-----------------------------------------------------------------------------------
\begin{thm}
Let $(F,G,E),(F_1,G_1,E),(F_2,G_2,E)\in\mathbb{BS}$. Then,
\begin{enumerate}[\it i.]
\item $(F,G,E)\sqcap\big[(F_1,G_1,E)\sqcup(F_2,G_2,E)\big]=
\big[(F,G,E)\sqcap(F_1,G_1,E)\big]\sqcup\big[(F,G,E)\sqcap(F_2,G_2,E)\big]$
\item $(F,G,E)\sqcup\big[(F_1,G_1,E)\sqcap(F_2,G_2,E)\big]=
\big[(F,G,E)\sqcup(F_1,G_1,E)\big]\sqcap\big[(F,G,E)\sqcup(F_2,G_2,E)\big]$
\end{enumerate}
\end{thm}
\proof It can be proved simply.

%-----------------------------------------------------------------------------------
\begin{thm}
Let $(F,G,E),(F_1,G_1,E),(F_2,G_2,E)\in\mathbb{BS}$. Then,
\begin{enumerate}[\it i.]
\item $\big((F,G,E)^{\tilde c}\big)^{\tilde c}=(F,G,E)$
\item $(\Phi,\tilde U,E)^{\tilde c}=(\tilde U,\Phi,E)$
\item $(\tilde U,\Phi,E)^{\tilde c}=(\Phi,\tilde U,E)$
\end{enumerate}
\end{thm}
\proof It is obvious from Definition \ref{d-complement}.

%-----------------------------------------------------------------------------------
\begin{thm}\label{t-demorganlaws}
Let $(F_1,G_1,E),(F_2,G_2,E)\in\mathbb{BS}$. Then, De Morgan's law is valid.
\begin{enumerate}[\it i.]
\item $\big[(F_1,G_1,E)\sqcup(F_2,G_2,E)\big]^{\tilde c}=
(F_1,G_1,E)^{\tilde c}\sqcap(F_2,G_2,E)^{\tilde c}$
\item $\big[(F_1,G_1,E)\sqcap(F_2,G_2,E)\big]^{\tilde c}=
(F_1,G_1,E)^{\tilde c}\sqcup(F_2,G_2,E)^{\tilde c}$
\end{enumerate}
\end{thm}
\pagebreak
\proof
\begin{enumerate}[\it i.]
\item Let $(F_1,G_1,E)\sqcup(F_2,G_2,E)=(H,K,E)$ and $(H,K,E)^{\tilde c}=(S,T,E)$.
Then, $H(e)=F_1(e)\cup F_2(e)$ and $K(f(e))=G_1(f(e))\cap G_2(f(e))$  for all
$e\in E_1$. Thus, we have
\begin{equation}\label{eq1}
S(e)=K(f(e))=G_1(f(e))\cap G_2(f(e))
\end{equation}
and
\begin{equation}\label{eq2}
T(f(e))=H(e)=F_1(e)\cup F_2(e)
\end{equation}
for all $e\in E_1$. Moreover, $(F_1,G_1,E)^{\tilde c}=(H_1,K_1,E)$ and
$(F_2,G_2,E)^{\tilde c}=(H_2,K_2,E)$. Then, $H_1(e)=G_1(f(e))$, $H_2(e)=G_2(f(e))$,
$K_1(f(e))=F_1(e)$ and $K_2(f(e))=F_2(e)$. Therefore,
\begin{equation}\label{eq3}
(H_1\sqcap H_2)(e)=G_1(f(e))\cap G_2(f(e))
\end{equation}
and
\begin{equation}\label{eq4}
(K_1\sqcap K_2)(e)=F_1(e)\cup F_2(e)
\end{equation}
for all $e\in E_1$. On the other hand, it can be seen clearly that right hand of
\eqref{eq1}  is equal to right hand of \eqref{eq3} and right hand of \eqref{eq2}
is equal to right hand of \eqref{eq4}. So, the proof is completed.
\item It can be proved similar way (\emph{i.}).
\end{enumerate}

%-----------------------------------------------------------------------------------
\begin{defn} \label{d-andproduct}
Let $(F_1,G_1,E),(F_2,G_2,E)\in\mathbb{BS}$. Then, \emph{and}-product of bipolar soft
sets $(F_1,G_1,E)$ and $(F_2,G_2,E)$ is, denoted by $(F_1,G_1,E)\wedge(F_2,G_2,E)$,
defined by $(F_1\wedge F_2)(e,e')=F_1(e)\cap F_2(e')$  and
$(G_1\wedge G_2)(f(e))=G_1(f(e))\cup G_2(f(e),f(e'))$ for all $e,e'\in E_1$.
\end{defn}

%-----------------------------------------------------------------------------------
\begin{defn} \label{d-orproduct}
Let $(F_1,G_1,E),(F_2,G_2,E)\in\mathbb{BS}$. Then, \emph{or}-product
of bipolar soft sets $(F_1,G_1,E)$ and $(F_2,G_2,E)$,  denoted by
$(F_1,G_1,E)\vee(F_2,G_2,E)$, is defined by $(F_1\vee
F_2)(e,e')=F_1(e)\cup F_2(e')$  and $(G_1\vee
G_2)(f(e),f(e'))=G_1(f(e))\cap G_2(f(e'))$ for all $e,e'\in E_1$.
\end{defn}

%-----------------------------------------------------------------------------------
\begin{exmp}Let $U=\{u_l,u_2,u_3,u_4,u_5,u_6,u_7,u_8\}$ be the universe which are eight houses
and $E=\{e_1,e_2,e_3,e_4,e_5,e_6\}$ be the set of parameters. Here,
$e_i$ $(i=1,2,3,4)$ stand for the parameters ``\emph{large}'',
``\emph{small}'', ``\emph{cheap}'',``\emph{modern}'',
``\emph{standard}'', ``\emph{expensive}'', respectively.  Therefore,
we can chose $E_1$ and $E_2$ sets as $E_1=\{e_1,e_3,e_5\}$ and
$E_2=\{e_2,e_4,e_6\}$. Now, we define the bijective function $f$ as
$f(e_i)=\neg e_i$ $(i=1,3,5)$. Here, the notion $\neg e_i$ means
that ``\emph{not} $e_i$'' for all $i=1,3,5$. Thus, we have following
results.
\begin{eqnarray*}
& & f(e_1) = \neg e_1=e_2 \\
& & f(e_3) = \neg e_3=e_4 \\
& & f(e_5) = \neg e_5=e_6 \\
\end{eqnarray*}
So, we can describe following the bipolar soft sets $(F_1,G_1,E)$,
$(F_2,G_2,E)$ and $(F_3,G_3,E)$.
\begin{eqnarray*}
(F_1,G_1,E) & = & \Big\{
\big\langle(e_1,\{u_1,u_3,u_4\}),(e_2,\{u_2,u_6\})\big\rangle,\\
&&\big\langle(e_3,\{u_2,u_5,u_7\}),(e_4,\{u_1,u_3,u_8\})\big\rangle,\\
&&\big\langle(e_5,\{u_3,u_4\}),(e_6,\{u_1,u_2,u_5,u_8\})\big\rangle\Big\}
\end{eqnarray*}
\begin{eqnarray*}
(F_2,G_2,E) & = & \Big\{
\big\langle(e_1,\{u_1,u_2,u_4\}),(e_2,\{u_3,u_5,u_6,u_7\})\big\rangle,\\
&&\big\langle(e_3,\{u_2,u_5\}),(e_4,\{u_1,u_3,u_4,u_8\})\big\rangle,\\
&&\big\langle(e_5,\{u_1,u_3,u_4\}),(e_6,\{u_2,u_5,u_7,u_8\})\big\rangle\Big\}.
\end{eqnarray*}
Then,

\begin{eqnarray*}
(F_1,G_1,E)\wedge(F_2,G_2,E)&=&\big\{\big\langle((e_1,e_1),\{u_1,u_4\}),((e_2,e_2),\{u_2,u_3,u_5,u_6,u_7\})\big\rangle,\\
&&\big\langle((e_1,e_3),\emptyset),((e_2,e_4),\{u_1,u_2,u_3,u_4,u_6,u_8\})\big\rangle,\\
&&\big\langle((e_1,e_5),\{u_1,u_3,u_4\}),((e_2,e_6),\{u_2,u_5,u_6,u_7,u_8\})\big\rangle,\\
&&\big\langle((e_3,e_1),\emptyset),((e_4,e_2),\{u_1,u_2,u_3,u_4,u_6,u_8\})\big\rangle,\\
&&\big\langle((e_3,e_3),\{u_2,u_5\}),((e_4,e_4),\{u_1,u_3,u_4,u_8\})\big\rangle,\\
&&\big\langle((e_3,e_5),\emptyset),((e_4,e_6),\{u_1,u_2,u_3,u_5,u_7,u_8\})\big\rangle,\\
&&\big\langle((e_5,e_1),\{u_1,u_3,u_4\}),((e_6,e_2),\{u_2,u_5,u_6,u_7,u_8\})\big\rangle,\\
&&\big\langle((e_5,e_3),\emptyset),((e_6,e_4),\{u_1,u_2,u_3,u_5,u_7,u_8\})\big\rangle,\\
&&\big\langle((e_5,e_5),\{u_3,u_4\}),((e_6,e_6),\{u_1,u_2,u_5,u_7,u_8\})\big\rangle\big\}
\end{eqnarray*}
and
\begin{eqnarray*}
(F_1,G_1,E)\vee(F_2,G_2,E)&=&\big\{\big\langle((e_1,e_1),\{u_1,u_2,u_3,u_4\}),((e_2,e_2),\{u_6\})\big\rangle,\\
&&\big\langle((e_1,e_3),\{u_1,u_2,u_3,u_4\}),((e_2,e_4),\emptyset\})\big\rangle,\\
&&\big\langle((e_1,e_5),\{u_1,u_3,u_4\}),((e_2,e_6),\{u_2\})\big\rangle,\\
&&\big\langle((e_3,e_1),\{u_1,u_2,u_4,u_5,u_7\}),((e_4,e_2),\{u_3\})\big\rangle,\\
&&\big\langle((e_3,e_3),\{u_2,u_5,u_7\}),((e_4,e_4),\{u_1,u_3,u_8\})\big\rangle,\\
&&\big\langle((e_3,e_5),\{u_1,u_2,u_3,u_4,u_5,u_7\}),((e_4,e_6),\{u_8\})\big\rangle,\\
&&\big\langle((e_5,e_1),\{u_1,u_3,u_4\}),((e_6,e_2),\{u_5\})\big\rangle,\\
&&\big\langle((e_5,e_3),\{u_2,u_3,u_4,u_5\}),((e_6,e_4),\{u_1,u_8\})\big\rangle,\\
&&\big\langle((e_5,e_5),\{u_1,u_3,u_4\}),((e_6,e_6),\{u_2,u_5,u_8\})\big\rangle\big\}
\end{eqnarray*}
\end{exmp}

%-----------------------------------------------------------------------------------
\begin{thm}
Let $(F_1,G_1,E),(F_2,G_2,E)\in\mathbb{BS}$. Then,
\begin{enumerate}[\it i.]
\item $\big[(F_1,G_1,E)\vee(F_2,G_2,E)\big]^{\tilde c}=
(F_1,G_1,E)^{\tilde c}\wedge(F_2,G_2,E)^{\tilde c}$
\item $\big[(F_1,G_1,E)\wedge(F_2,G_2,E)\big]^{\tilde c}=
(F_1,G_1,E)^{\tilde c}\vee(F_2,G_2,E)^{\tilde c}$
\end{enumerate}
\end{thm}
\proof The proof is clear.
%\begin{enumerate}[\it i.]
%\item Let $E=E_1\cup E_2$, $(F_1,G_1,E)\vee(F_2,G_2,E)=(H,K,E_1\times E_1)$ and $(H,K,E_1\times E_1)^{\tilde c}=(S,T,E_1\times E_1)$.
%Then, $H(e,e')=F_1(e)\cup F_2(e')$ and $K(f(e))=G_1(f(e))\cap G_2(f(e'))$  for all
%$(e,e')\in E_1\times E_1$. Thus, we have
%\begin{equation}\label{eq1}
%S(e)=K(f(e))=G_1(f(e))\cap G_2(f(e'))
%\end{equation}
%and
%\begin{equation}\label{eq2}
%T(f(e))=H(e)=F_1(e)\cup F_2(e')
%\end{equation}
%for all $(e,e')\in E_1\times E_1$. Moreover, $(F_1,G_1,E)^{\tilde c}=(H_1,K_1,E)$ and
%$(F_2,G_2,E)^{\tilde c}=(H_2,K_2,E)$. Then, $H_1(e)=G_1(f(e))$, $H_2(e)=G_2(f(e))$,
%$K_1(f(e))=F_1(e)$ and $K_2(f(e))=F_2(e)$. Therefore,
%\begin{equation}\label{eq3}
%(H_1\sqcap H_2)(e)=G_1(f(e))\cap G_2(f(e))
%\end{equation}
%and
%\begin{equation}\label{eq4}
%(K_1\sqcap K_2)(e)=F_1(e)\cup F_2(e)
%\end{equation}
%for all $e\in E_1$. On the other hand, it can be seen clearly that right hand of
%\eqref{eq1}  is equal to right hand of \eqref{eq3} and right hand of \eqref{eq2}
%is equal to right hand of \eqref{eq4}. So, the proof is completed.
%\item It can be proved similar way (\emph{i.}).
%\end{enumerate}

%***********************************************************************************
%% section: An applications of bipolar soft sets in decision making
\section{An applications of bipolar soft sets in decision making}
In this section we will construct a decision making  method over the
bipolar soft set. Firstly, we will define  some notions that necessary
 to construct algorithm of decision making method.
%------------------------------------------------------------
\begin{defn} Let $E=\{e_1,e_2,...,e_n\}$ be a parameter set,
$U=\{u_1,u_2,...,u_m\}$ be initial universe and $(F,G,E)$ be a BSS over $U$.Then,
score of an object, denoted by $s_i$, is
computed as $s_i = c_i^+-c_i^- $.
Here, $c_i^+$ and $c_i^-$  are computed form $c_i^+=\sum_{j=1}^na_{ij}$ and
$c_i^-=\sum_{j=1}^nb_{ij}$
\end{defn}

 Now we present an algorithm for
most appropriate selection of an object.

\section*{Algorithm}

\textbf{Step 1:} Input the bipolar soft Set $(F,G,E)$\\
\textbf{Step 2:} Consider the bipolar soft set $(F,G,E)$ and write it in tabular form\\
\textbf{Step 3:} Compute the score $s_i$ of $h_i$ $\forall i$\\
\textbf{Step 4:} Find  $s_k = maxs_i$\\
\textbf{Step 5:} If k has more than one value then any one
of $h_i$ could be the preferable choice.\\

Let us use the algorithm to solve the problem.

\begin{exmp}
We consider the problem to select the most suitable house which Mr.
$X$ is going to choose on the basis of his $10$ number of parameters out
of $m$ number of houses.
Let
$E_1=\{e_1=beautiful,e_2=cheap, e_3=in\, good \,repairing,\\
e_4=moderate,e_5=wooden\}$ and $E_2=\{ e_6=not \,beautiful,e_7=not\, cheap,\\
e_8=not\,in\, good\,repairing,e_9=not\, moderate,e_{10}=not\,wooden\}$ and $U$ be a initial universe. Here, we consider
$E=E_1\cup E_2$. For the sake of shortness, we will denote the
$(F(e_i),G(f(e_i)))$ with $FG(e_i)$

\textbf{Step 1:} Let us consider BSS $(F,G,E)$ over $U$ defined as follow,
\small
\begin{eqnarray*}
(F,G,E) & = & \Big\{
\big\langle(e_1,\{u_1,u_2\}),(e_5,\{u_3,u_5,u_6,u_7\})\big\rangle,\big\langle(e_2,\{u_2\}),(e_6,\{u_1,u_3,u_4,u_8\})\big\rangle,\\
&&\big\langle(e_3,\{u_1,u_3,u_4\}),(e_7,\{u_2,u_5,u_7,u_8\})\big\rangle,\big\langle(e_4,\{u_5\}),(e_8,\{u_2,u_3,u_4\})\big\rangle\\
&&\big\langle(e_5,\{u_1,u_3,u_8\}),(e_{10},\{u_5,u_6\})\big\rangle\Big\}\\
\end{eqnarray*}
\normalsize

\textbf{Step 2:}
Tabular representation of the BSS $(F,G,E )$ is as
below:
\footnotesize
\begin{center}
\begin{tabular}{c|ccccc}
$(F,G,E)$ & $FG(e_1)$ & $FG(e_2)$ & $FG(e_3)$ & $FG(e_4)$ & $FG(e_5)$\\
  \hline
$u_1$ & $(1,0)$ & $(0,1)$ & $(1,0)$ & $(0,0)$ & $(1,0)$\\
$u_2$ & $(1,0)$ & $(1,0)$ & $(0,1)$ & $(0,1)$ & $(0,0)$\\
$u_3$ & $(0,1)$ & $(0,1)$ & $(1,0)$ & $(0,1)$& $(1,0)$ \\
$u_4$ & $(1,0)$ & $(0,1)$ & $(1,0)$ & $(0,1)$ & $(0,0)$ \\
$u_5$ & $(0,1)$ & $(0,0)$ & $(0,1)$ & $(1,0)$ & $(0,1)$ \\
$u_6$ & $(0,1)$ & $(0,0)$ & $(0,0)$ & $(0,0)$ & $(0,1)$ \\
$u_7$ & $(0,1)$ & $(0,0)$ & $(0,1)$ & $(0,0)$& $(0,0)$\\
$u_8$ & $(0,0)$ & $(0,1)$ & $(0,1)$ & $(0,0)$ & $(1,0)$\\
\end{tabular}
\end{center}
\begin{center}
Tabular form of the BSS $(F,G,E)$
\end{center}
\normalsize

\textbf{Step 3:} Score of each elements of $U$ is obtained as follow
\begin{center}
\begin{tabular}{c|c|c|c}
    & $c_i^+$ & $c_i^-$ & $s_i$ \\
  \hline
  $u_1$ & 3 & 1 & 2 \\
  $u_2$ & 2 & 2 & 0 \\
  $u_3$ & 2 & 3 & -1 \\
  $u_4$ & 2 & 2 & 0 \\
  $u_5$ & 1 & 3 & -2  \\
  $u_6$ & 0 & 2 & -2 \\
  $u_7$ & 0 & 2 & -2 \\
  $u_8$ & 1 & 2 & -1 \\
\end{tabular}
\end{center}
\end{exmp}

\textbf{Step 4:} Since $maks_i=2$, optimal element of $U$ is $u_1$.

%***********************************************************************************
%% section: Conclusion
%\section{Conclusion}

%***********************************************************************************
%% bibliography

\end{document}